\documentclass[9pt,conference]{ieeeconf}

\IEEEoverridecommandlockouts                              

\usepackage{color}
\usepackage{probsoln,nicefrac}
\usepackage{graphics} 
\usepackage{epsfig} 
\usepackage{subfigure}
\usepackage{soul}
\usepackage{times} 
\usepackage{amsmath} 
\usepackage{amssymb}  
\usepackage[english]{babel}
\usepackage{mathtools}

\usepackage{mathrsfs}
\usepackage{enumerate}
\usepackage{tikz-cd}
\usepackage{cite}
\usepackage{float}

\newtheorem{theorem}{Theorem}[section]
\newtheorem{proposition}{Proposition}[section]
\newtheorem{definition}{Definition}[section]
\newtheorem{lemma}[definition]{Lemma}
\newtheorem{example}[definition]{Example}

\newtheorem{remarkth}[definition]{Remark}
\newenvironment{remark}{\begin{remarkth}\upshape}{\end{remarkth}}

\newcommand{\R}{\mathbb{R}}

\newcommand{\F}{\mathbb{F}}

\newcounter{met}

\newcommand{\norm}[1]{\left\lVert#1\right\rVert}

\author{Mar\'ia Emma Eyrea Iraz\'u, Asier L\'opez-Gord\'on, Manuel de Le\'on, Leonardo J. Colombo
\thanks{M.~E.~E.~I.~(maemma@mate.unlp.edu.ar) is with CONICET and CMaLP, Dept. Mathematics, UNLP, C/ 1 y 115, La Plata 1900, Buenos Aires, Argentina.  L.~C.~(leonardo.colombo@car.upm-csic.es) is with CAR (CSIC-UPM), Ctra. M300 Campo Real, Km 0, 200, Arganda del Rey - 28500 Madrid, Spain. M.~d.~L.~(mdeleon@icmat.es) and A.~L.-G.~(asier.lopez@icmat.es) are with ICMAT (CSIC-UAM-UC3M-UCM), C/ Nicol\'as Cabrera, 15 - 28049 Madrid, Spain.}
}%

\title{\LARGE \bf
Hybrid Routhian reduction for \\ simple hybrid forced Lagrangian systems
}

\begin{document}

\maketitle
\thispagestyle{empty}
\pagestyle{empty}


\begin{abstract}
This paper discusses Routh reduction for simple hybrid forced mechanical systems. We give general conditions on whether it is possible to perform symmetry reduction for a simple hybrid Lagrangian system subject to non-conservative external forces, emphasizing the case of cyclic coordinates. We illustrate the applicability of the symmetry reduction procedure with an example and numerical simulations.

\end{abstract}

\section{Introduction}

Dimensionality reduction for large scale systems has become an active problem of interest within the automatic control and robotics communities. For instance, in large robotic swarms, guidance and trajectory planning algorithms for coordination while optimizing qualitative features for the swarm of multiple robots are determined by solutions of nonlinear equations which demand a high-computational cost along its integration. A key element in the reduction is a Lie group of symmetries. Lie groups of symmetries appear naturally in many control systems problems \cite{Bl}, \cite{grizzle2}. Examples of invariant control problems on Lie groups include motion planning for underwater vehicles \cite{Leonard1}, conflict resolution in differential games \cite{tomlin}, collective motion in biological models \cite{JK2}, and coordination of multi-agent systems \cite{bonnabel}, \cite{leo1}. The construction of methods for reduction of dimensionality also permits, for instance, fast computations for the generation of optimal trajectories in optimal control problems of mechanical systems \cite{leo2}.

Hybrid systems are dynamical systems with continuous-time and discrete-time components on its dynamics. Simple hybrid systems are a class of hybrid system introduced in \cite{SHS}, denoted as such because of their simple structure. A simple hybrid system is characterized by a tuple $\mathscr{L}=(D, X, S, \Delta)$ where $D$ is a smooth manifold, $X$ is a smooth vector field on $D$, $S$ is an embedded submanifold of $D$ with co-dimension $1$, and $\Delta:S\to D$ is a smooth embedding. This type of hybrid system has been mainly employed for the understanding of locomotion gaits in bipeds and insects \cite{ames2006}, \cite{HoFuKoGu}, \cite{Biped-book}. In the situation where the vector field $X$ is associated with a mechanical system (Lagrangian or Hamiltonian), alternative approaches for mechanical systems with nonholonomic and unilateral constraints have been considered in \cite{will}, \cite{cortes}, \cite{cortes2}, \cite{ibort}, \cite{ibort2}.

When a dynamical system exhibits a symmetry, it produces a conserved quantity for the system and allows to reduce the degrees of freedom in the dynamics by using the conserved quantity. One of the classical reduction by symmetry procedures in mechanics is the Routh reduction method \cite{goldstein}, \cite{Foundations}. During the last few years there has been a growing interest in Routh reduction, mainly motivated by physical applications - see \cite{Edu}, \cite{quasi}, \cite{Routhstages} and references therein. Routh reduction for hybrid systems has been introduced and applied in the field of bipedal locomotion \cite{ames2006}. The reduced simple hybrid system is called simple hybrid Routhian system \cite{amesrouth}. A hybrid scheme for Routh reduction for simple hybrid Lagrangian systems with cyclic variables is found in~\cite{amesrouth} and \cite{NAHS}, inspired to gain a better understanding of bipedal walking models (see also~\cite{ames2006} and references therein). Symplectic reduction for hybrid Hamiltonian systems has been introduced in \cite{amesham} and extended Poisson reduction in \cite{emmapoisson} and to time-dependent systems in \cite{tdep}, but to the best of our knowledge, the hybrid analogue for symmetry reduction in mechanical systems subject to external forces has not been explored in the literature. This paper attempts to go one step further and to consider symmetry reduction of simple hybrid Lagrangian system subject to external non-conservative forces via Routh reduction for simple hybrid forced systems. Fundamental to the reduction procedure has been the recent work \cite{manololainzlopen} on reduction of (non-hybrid) forced Lagrangian systems. 

The paper is organized as follows. Sec. II presents
the necessary background on the geometry of forced mechanical systems and Routh reduction. Sec. III introduces the class of simple hybrid forced Lagrangian and Hamiltonian systems under consideration and the corresponding relation between both formalisms. The reduction scheme is proposed in Section IV. The reduction technique has been illustrated both analytically and numerically in an expository example.

\section{Routh reduction for forced systems}
We start by recalling some basic facts about mechanical systems subject to external forces.

Let $Q$ be an $n$-dimensional differentiable
manifold with local coordinates $(q^i)$, $1\leq i\leq n$, the
configuration space of a mechanical system. Denote by $TQ$ its
tangent bundle, that is, if $T_{q}Q$ denotes the tangent space of $Q$ at the point $q\in Q$, then $\displaystyle{TQ\coloneqq \bigcup_{q\in Q}T_{q}Q}$, with induced local coordinates $(q^i, \dot{q}^i)$. Since $T_{q}Q$ has a vector space structure, we may consider its dual space, $T^{*}_{q}Q,$ and define the cotangent bundle as $\displaystyle{T^{*}Q\coloneqq \bigcup_{q\in Q}T^{*}_{q}Q},$ with local coordinates $(q^i,p_i)$.

The dynamics of a mechanical system can be determined by the Euler-Lagrange equations associated with a Lagrangian function $L:TQ\to\mathbb{R}$ given by $L(q,\dot{q})=K(q,\dot{q})-V(q),$ where $K=\frac{1}{2} \norm{\dot q}_q^2$ is the kinetic energy and $V:Q\to\mathbb{R}$ the potential energy. The Lagrangian $L$ is said to be regular if $\displaystyle{\det\mathcal{M}\coloneqq \det\left(\frac{\partial^2L}{\partial\dot{q}^{i}\partial\dot{q}^{j}}\right)\neq0}$ for all $i,j$ with $1\leq i,j\leq n$.

The equations describing the dynamics of the system are given by the Euler-Lagrange equations $\displaystyle{\frac{d}{dt}\left(\frac{\partial L}{\partial\dot{q}^{i}}\right)=\frac{\partial L}{\partial q^{i}}}$, with $i=1,\ldots,n;$ a system of $n$ second-order ordinary differential equations. If $L$ is regular, the Euler-Lagrange equations induce a vector field $X_L:TQ\to T(TQ)$ describing the dynamics of the Lagrangian system, given by $$X_L(q^i,\dot{q}^i)=\left(q^i,\dot{q}^i;\dot{q}^i, \mathcal{M}^{-1}\left(\frac{\partial L}{\partial q^i}-\frac{\partial^2L}{\partial\dot{q}^{i}\partial q^j}\dot{q}^{j}\right)\right).$$

In addition, the motion of the system  may be influenced by a non-conservative force (conservative forces may be included into the potential energy), which is a smooth map $F:TQ\to T^{*}Q,$ locally given by $F=F_{i}dq^i$ and geometrically respresented by a $1$-form on $Q$.  At a given position and velocity, the force will act against variations of the position (i.e., virtual displacements). Lagrange-d'Alembert principle leads to the so-called \textit{forced Euler-Lagrange equations} $\displaystyle{\frac{d}{dt}\left(\frac{\partial L}{\partial\dot{q}^i}\right)-\frac{\partial L}{\partial q^i}=F_{i}(q,\dot{q})}$, $i=1,\ldots,n$; a system of $n$ second-order ordinary differential equations. If $L$ is regular equations those equations induce a vector field $X_{L}^{F}:TQ\to T(TQ)$ describing the dynamics of the forced Lagrangian system, given by
$$X_{L}^{F}(q^i,\dot{q}^i)=\left(q^i,\dot{q}^i;\dot{q}^i, \mathcal{M}^{-1}\left(F_{i}+\frac{\partial L}{\partial q^i}-\frac{\partial^2L}{\partial\dot{q}^{i}\partial q^j}\dot{q}^{j}\right)\right).$$

For the Lagrangian $L:TQ\to\mathbb{R}$, let us denote by $\F L\colon TQ\to  T^*Q$ the Legendre transformation associated with $L$; that is, $\F L\colon TQ\to  T^*Q,\, (q,\dot q)\mapsto (q,p\coloneqq \partial L/\partial \dot q)$. The map $\F L:TQ\to T^*Q$ relates velocities and momenta. In fact, the Legendre Transformation connects Lagrangian and Hamiltonian formulations of mechanics (see \cite{Foundations}).

We said that the Lagrangian is \textit{hyperregular}, if $\F L$ is a diffeomorphism between $TQ$ and $T^*Q$ (this is always the case for mechanical Lagrangians). If $L$ is hyperregular, one can work out the velocities $\dot q=\dot q(q,p)$ in terms of $(q,p)$ and define the Hamiltonian function (the ``total energy'') $H\colon T^*Q\to \R$ as
$H(q,p)=p^T \dot q(q,p)- L(q,\dot q(q,p))$, where we have used the inverse of the Legendre transformation to express $\dot q=\dot q(q,p)$. The evolution vector field corresponding to the Hamiltonian $H$, denoted by $X_{H}$, is defined by 
$\displaystyle{
X_H=\frac{\partial H}{\partial p_i}\frac{\partial}{\partial q^i}-\frac{\partial H}{\partial q^i}\frac{\partial}{\partial p_i}}$,
and its integral curves are solutions of Hamilton's equations $\displaystyle{
\dot q^i=\frac{\partial H}{\partial p_i},\,\dot p_i= -\frac{\partial H}{\partial q^i}}$.

If the system is influenced by a nonconservative force, \textit{forced Hamilton's equations} are given by
$\displaystyle{
\dot q^i=\frac{\partial H}{\partial p_i},\quad \dot p_i= F^{H}_{i}-\frac{\partial H}{\partial q^i}}$, where $F^{H}=\F L(F)$.

There exists a large class of systems for which the Lagrangian (resp. Hamiltonian) does not depend on some of the generalized coordinates. Such coordinates are called \textit{cyclic}  or \textit{ignorable}, and the corresponding generalized momenta are easily checked to be constants of motion - see \cite{Foundations}, \cite{goldstein}. Routh's reduction procedure is a classical reduction technique which takes advantage of the conservation laws to define a \textit{reduced Lagrangian} function, so-called \textit{Routhian funtion}, such that the solutions of the Euler-Lagrange equations for the Routhian are in correspondence with the solutions of Euler-Lagrange equations for the original Lagrangian, when the conservation of momenta is taken into account. 

Routh's reduction can be extended to forced systems as follows \cite{goldstein}:
the starting point for Routhian reduction for forced Lagrangian systems is a configuration space of the form $Q=Q_{1}\times Q_{2}$, where we denote an element $q^{i}\in Q$ by $q^{i}=(q^{1},q^j)$ with $q^1\in Q_1$ and $q^j\in Q_{2}$, with $j=2,\ldots,n$.

Let $L(q^1,\dot q^1,q^j,\dot q^j)$ be a hyper-regular Lagrangian with cyclic coordinate $q^1$, that is, $\displaystyle{\frac{\partial L}{\partial q^1}=0}$ and let $F_{i}$ be a non-conservative force such that $F_{i}$ is independent of $q^1$ for all $i=1,...,n$ and $F_{1}(q^2,...,q^n)=0$. Fundamental to reduction is the notion of a momentum map $J_L:TQ\rightarrow \mathfrak{g}^{*}$, which makes explicit the conserved quantities in the system, with $\mathfrak{g}$ the Lie algebra associated with the Lie group of symmetries and $\mathfrak{g}^{*}$ denoting its dual as vector space. In the framework we are considering here, $J_L(q^1,\dot q^1,q^j,\dot q^j)=\displaystyle{\frac{\partial L}{\partial \dot{q}^1}}.$ Fix a value of the momentum $\mu=\displaystyle{\frac{\partial L}{\partial \dot{q}^1}}$. Since $L$ is  hyper-regular, the last equation admits an inverse, and allows us to write $\dot{q}^{1}=f(q^2,...,q^n,\dot{q}^2,...,\dot{q}^n,\mu)$. Consider the function $R^{\mu}(q^j,\dot q^j)=\left(L-\dot q^1\mu\right)\Big{|}_{\mu}$, where the notation $\mid_{\mu}$ means that we have used the relation $\mu=\displaystyle{\frac{\partial L}{\partial \dot{q}^1}}$ to replace all the appearances of $\dot q^1$ in terms of $(q^j,\dot q^j)$ and the parameter $\mu$. $R^{\mu}$ is called \textit{Routhian function}.

 If we regard $R^{\mu}$ and $F_{\mu}$ as a new Lagrangian and external force in the variables $(q^j,\dot q^j)$, then the solutions of the forced Euler-Lagrange equations for $R^{\mu}$ with $F_{\mu}$ are in correspondence with those of $L$ with $F$ when one takes into account the relation $\mu=\displaystyle{\frac{\partial L}{\partial \dot{q}^1}}$. More precisely:

	\begin{itemize}
		\item[(a)] Any solution of the forced Euler-Lagrange equation for $L$ and $F$ with momentum $\mu=\displaystyle{\frac{\partial L}{\partial \dot{q}^1}}$ projects onto a solutions of the forced Euler-Lagrange equations for $R^{\mu}$ and $F_{\mu}$, $\frac{d}{dt}\left(\frac{\partial R^{\mu}}{\partial\dot{q}^{j}}\right)-\frac{\partial R^{\mu}}{\partial q^{j}}=(F_{\mu})_j$, $j=2,\ldots,n$. These equations will be referred to as \textit{forced Routh equations} and they induce a vector field $X_{R}^{F}:TQ_{2}\to T(TQ_{2})$ describing the dynamics of the reduced system, called \textit{Routhian vector field}.
		\item[(b)] Conversely, any solution of forced Routh equations for $R^{\mu}$ and $F_{\mu}$ can be lifted to a solution of the forced Euler-Lagrange equations for $L$ and $F$ with $\mu=\displaystyle{\frac{\partial L}{\partial \dot{q}^1}}$.
	\end{itemize}

\begin{example}[Billiard with dissipation]\label{example}

Consider a particle of mass $m$ in the plane which is free to move inside the surface defined by $x^2+y^2=1$. The surface of the ``billiard'' is assumed to be rough in such a way that the friction is non-linear on the velocities.

The Lagrangian $L:T\mathbb{R}^{2}\to\mathbb{R}$ is given by 
$\displaystyle{
L(x,y,\dot{x},\dot{y})=\left[\frac{1}{2}m(\dot x^2+\dot y^2)\right]}$and $F(x,y,\dot{x},\dot{y})=F_{x}dx+F_{y}dy$ is an external force given by $F_{x}=2c(\dot{x}xy-\dot{y}x^2)$, $F_{y}=-2c(\dot{y}xy-\dot{x}y^2)$, for a constant $c>0$. The equations of motion for the particle off the boundary are then 
$$
m\ddot x=-2c(\dot{y}x^2-\dot{x}xy), \,  m\ddot y=2c(\dot{x}y^2-\dot{y}xy).$$

By introducing polar coordinates $L$ and $F$ become
$\displaystyle{
L(\theta,r,\dot{\theta},\dot{r})=\frac{m}{2}(\dot r^2+ r^2\dot \theta^2)}$, and $F(\theta,r,\dot{\theta},\dot{r})=-2cr^3\dot{\theta}dr+0d\theta$, respectively. $L$ is hyper-regular and both $L$ and $F$ are independent of $\theta$. The forced Euler-Lagrange equations (in polar coordinates) are $\displaystyle{\ddot r=(r-\frac{2c}{m}r^3)\dot{\theta},\,\,mr^2\ddot{\theta}=0}$.

Note that the momentum map $J_L$ for $\theta$, $J_L(r,\dot r,\theta,\dot \theta)= mr^2\dot\theta$ is preserved, that is, by considering $\mu=mr^2\dot{\theta}$ (i.e., $\dot{\theta}=\frac{\mu}{mr^2}$) the Routhian and the reduced force takes the form $R^{\mu}(r,\dot{r})=\frac{m}{2}\dot r^2-\frac{\mu^2}{2mr^2}, \,\, F_{\mu}=-2cr\frac{\mu} {m}dr$ and the forced reduced Euler-Lagrange equations for the Routhian $R^{\mu}$ are given by $\ddot r=\frac{\mu^2}{m^2r^3}-2cr\frac{\mu} {m^2}$.
\end{example}

\section{Simple hybrid forced Lagrangian systems}\label{Sec3.A}
Roughly speaking, the term \textit{hybrid system} refers to a dynamical system which exhibits both a continuous and a discrete time behaviors. In the literature, one finds slightly different definitions of hybrid system depending on the specific class of applications of interest. For simplicity, and following \cite{SHS}, ~\cite{amesrouth}, we will restrict ourselves to the so-called simple hybrid mechanical systems in Lagrangian form. 

Simple hybrid systems \cite{SHS} (see also \cite{Biped-book}) are characterized by the 4-tuple $\mathscr{L}=(D, X, S, \Delta)$ where $D$ is a smooth manifold, the \textit{domain}, $X$ is a smooth \textit{vector field} on $D$, $S$ is an embedded submanifold of $D$ with co-dimension $1$ called \textit{switching surface}, and $\Delta:S\to D$ is a smooth embedding called the \textit{impact map}. $S$ and $\Delta$ are also refered as the \textit{guard} and \textit{reset map}, respectively. The dynamics associated with a simple hybrid system is described by an autonomous system with impulse effects as in \cite{Biped-book}. We denote by $\Sigma_{\mathscr{L}}$ the \textit{simple hybrid dynamical system} generated by $\mathscr{L}$, that is,  \begin{equation}\label{LHS}\Sigma_{\mathscr{L}}:\begin{cases} \dot{\gamma}(t)=X(\gamma(t)),\quad\quad \gamma^{-}(t)\notin{S} \\ \gamma^{+}(t)=\Delta(\gamma^{-}(t))\quad \gamma^-(t)\in{S} \end{cases}\end{equation} where $\gamma:I\subset\mathbb{R}\to D$, and $\gamma^{-}$, $\gamma^{+}$ denote the states immediately before and after the times when integral curves of $X$ intersects ${S}$ (i.e., pre and post impact of the solution $\gamma(t)$ with ${S}$), where $\gamma^{-}(t)\coloneqq \displaystyle{\lim_{\tau\to t^{-}}}x(\tau)$,\, $\gamma^{+}(t)\coloneqq \displaystyle{\lim_{\tau\to t^{+}}}x(\tau)$ are the left and right limits of the state trajectory $\gamma(t)$.

A solution of a hybrid system may experience a Zeno state if infinity many impacts occur in a finite amount of time. It is particularly problematic in applications where numerical work is used, as computation time grows infinitely large at these Zeno points.  There are two primary modes through which Zeno behavior can occur: (i) A trajectory is reset back onto the guard, prompting additional resets. To exclude this behavior, we require that $S \cap \overline{\Delta}(S) = \emptyset$, where $\overline{\Delta}(S)$ denotes the closure as a set of $\Delta(S)$. This ensures that the trajectory will always be reset to a point with positive distance from the guard. (ii) The set of times where a solution to our system reaches the guard (and is correspondingly reset) has a limit point. This happens, for example, in the case of the \textit{bouncing ball} with coefficient of restitution $1/2$ - see \cite{brogliato}. To exclude these types of situations, we require the set of impact times to be closed and discrete, as in \cite{Biped-book}, so, we will assume implicitly throughout the remainder of the paper that   $\overline{\Delta}({S})\cap{S}=\emptyset$ and the set of impact times is closed and discrete. 

Given a smooth (constraint) function $h:Q\rightarrow \R$ on a configuration space $Q$ such that $h^{-1}(0)$ is a smooth submanifold, we can construct a domain and a guard explicitly - see \cite{brogliato}, \cite{amesrouth}. For this constraint function we have an associated domain, $D$,
defined to be the manifold (with boundary) $D=\{(q,\dot{q})\in TQ:h(q)\geq 0\}$. Similarly, we have an associated guard, $S$, defined as the submanifold of $D$ as
$S=\{(q,\dot{q})\in TQ:h(q)=0,\, dh_{q}\dot{q}\geq 0\}$, where $dh_{q}=\frac{\partial h}{\partial q}$. 

In a simple hybrid Lagrangian system the impact can be obtained from the  Newtonian impact equation (see \cite{brogliato} for instance) $P:TQ\rightarrow TQ$ given by
$$P(q,\dot{q})=\dot{q}-(1+e)\frac{dh_{q}\dot{q}}{dh_{q}M(q)^{-1}dh_{q}^T}M(q)^{-1}dh_{q}^T$$
where $M(q)$ is the inertial matrix for the Lagrangian system.

\begin{definition}
	A simple hybrid system $\mathscr{L}=(D, X, {S}, \Delta)$ is said to be a \textit{simple hybrid forced Lagrangian system} if it is determined by $\mathscr{L}_{F}\coloneqq (TQ, X_{L}^{F}, {S}, \Delta)$, where $X_{L}^{F}:TQ\to T(TQ)$ is the forced Lagrangian vector field, ${S}$ the switching surface and $\Delta:{S}\to TQ$ the impact map.\end{definition}

\begin{definition} The \textit{simple hybrid forced dynamical system} generated by $\mathscr{L}_{F}$ is given by
	\begin{equation}\label{RHDS}\Sigma_{\mathscr{L}_{F}}:\begin{cases} \dot{\gamma}(t)=X_{L}^{F}(\gamma(t)), \hbox{ if } \gamma^{-}(t)\notin{S},\\ \gamma^{+}(t)=\Delta(\gamma^{-}(t)),\hbox{ if } \gamma^-(t)\in{S}, \end{cases}\end{equation}where $\gamma(t)=(q^{a}(t),\dot{q}^{a}(t))\in TQ $.\end{definition}

\begin{definition}\label{def:flow} A \textit{hybrid flow} for $\mathscr{L}_{F}$ is a tuple $\chi^{\mathscr{L}_{F}}=(\Lambda,\mathcal{J},\mathscr{C})$, where
	\begin{itemize}
		\item $\Lambda=\{0,1,2,...\}\subseteq \mathbb{N}$ is a finite (or infinite) indexing set,
		\item $\mathcal{J}=\{I_{i}\}_{i\in \Lambda}$ a set of intervales, called hybrid intervals where $I_{i}=[\tau_{i},\tau_{i+1}]$ if $i, i+1\in \Lambda$ and $I_{N-1}=[\tau_{N-1},\tau_{N}]$ or $[\tau_{N-1},\tau_{N})$ or $[\tau_{N-1},\infty)$ if $|\Lambda|=N$, $N$ finite, with $\tau_{i},\tau_{i+1},\tau_{N}\in \R$ and $\tau_{i}\leq \tau_{i+1}$,
		\item $\mathscr{C}=\{c_{i}\}_{i\in \Lambda}$ is a collection of solutions for the vector field $X^{F}_L$ specifying the continous-time dynamics, i.e., $\dot{c_{i}}=X^{F}_L(c_{i}(t))$ for all $i\in \Lambda$,
		and such that for each $i,i+1\in \Lambda$, (i) $c_{i}(\tau_{i+1})\in S$, and (ii) $\Delta(c_{i}(\tau_{i+1}))=c_{i+1}(\tau_{i+1})$.
	\end{itemize}
\end{definition}

Analogously, one can introduce the notion of hybrid flow $\chi^{\mathscr{H}}$ for a simple hybrid forced Hamiltonian system $\mathscr{H}$. The relation between both hybrid flows is given by the following result, based on the well-known equivalence between the Lagrangian and Hamiltonian dynamics in the hyperregular case achieved via $\mathbb{F}L$.

\begin{proposition}\label{Proposition2}
If $\chi^{\mathscr{L}_{F}}=(\Lambda,\mathcal{J},\mathscr{C})$ is a hybrid flow for $\mathscr{L}_{F}$, $S_H=\F L(S)$, and $\Delta_H$ is defined in such a way that $\F L\circ \Delta=\Delta_H\circ \F L\mid_{S}$, then $\chi^{\mathscr{H}}=(\Lambda,\mathcal{J},(\F L)(\mathscr{C}))$ with $(\F L)(\mathscr{C})=\{(\F L)(c_{i})\}_{i\in \Lambda}$.
\end{proposition}

\textit{Proof:} If $c_i(t)$ is an integral curve of $X_{L}^{F}$, $\tilde c_i(t)=(\F L\circ c_i)(t)$ is an integral curve for $X_{H}^{F}$. In this way, if we consider a solution $c_0(t)$  with initial value $c_0=(q_0,\dot q_0)$ defined on $[\tau_0,\tau_1]$, then $\tilde c_0(t)$ is a solution with initial value $\tilde c_0=(q_0,p_0)$ defined on $[\tau_0,\tau_1]$. Likewise for a solution $c_1(t)$ defined on $[\tau_1,\tau_2]$, we get a corresponding solution $\tilde c_1(t)$ defined on the same hybrid interval $[\tau_1,\tau_2]$. Proceeding inductively, one finds $c_i(t)$ defined on  $[\tau_i,\tau_{i+1}]$. 	It only remains to check that $\tilde c_i(t)$ satisfies $\tilde c_i(\tau_{i+1})\in S_{H}$ and $\Delta_{H}(\tilde c_{i}(\tau_{i+1}))=\tilde c_{i+1}(\tau_{i+1})$, but using the properties of $\F L$,
\begin{itemize}
		\item[(i)] $\tilde c_i(\tau_{i+1})=(\F L\circ c_{i})(\tau_{i+1})=\F L(c_i(\tau_{i+1}))$ and given that $c_i(\tau_{i+1})\in S$ then $\tilde c_i(\tau_{i+1})\in S_{H}.$
		\item[(ii)]  $\Delta_H(\tilde c_{i}(\tau_{i+1}))=\Delta_H\circ \F L\circ c_{i}(\tau_{i+1})=\F L\circ \Delta\circ c_{i}(\tau_{i+1})=\F L\circ c_{i+1}(\tau_{i+1})=\tilde c_{i+1}(\tau_{i+1})$.\hfill$\square$
\end{itemize}

\begin{example}\label{ex2}
Continuing with Example \ref{example}, we consider the guard as the subset of $T\R^2\simeq\mathbb{R}^{2}\times\mathbb{R}^{2}$ given by
\[
S=(T\mathbb{R}^2)\cap \{x^2+y^2=1, (\dot x,\dot y)\cdot (x,y)\geq 0 \}. 
\]
 Under the assumption of an elastic collision, using the Newtonian impact equation, the reset map $(x,y,\dot x^{-},\dot y^{-})\mapsto (x,y,\dot x^{+},\dot y^{+})$, is given by
	\begin{align}
	\dot x^{+}&=\dot x^{-}  {-2 (x\dot x^{-}+y\dot y^{-})} x,\label{eq11}\\
	\dot y^{+}&=\dot y^{-} {-2 (x\dot x^{-}+y\dot y^{-})} y.\label{eq12}
	\end{align}
	Therefore, the 4-tuple $\mathscr{L}_{F}=(TQ,X_L^{F},S,\Delta)$, is a simple hybrid forced Lagrangian system with $Q=\R^{2}$, and $L$ and $F$ as described in Example \ref{example}.
\end{example}

\section{Routh reduction of simple hybrid forced Lagrangian systems}
Let $\mathscr{L}_{F}=(TQ,X_{L}^{F},S,\Delta)$ be a simple hybrid forced Lagrangian system. The starting point for symmetry reduction is a Lie group action $\psi\colon G\times Q\to Q$ of some Lie group $G$ on the manifold $Q$. We will assume that all the actions satisfy some regularity conditions as to do reduction, for instance, one can consider free and proper actions \cite{Foundations}.

There is a natural lift $\psi^{T^{*}Q}$ of the action $\psi$ to $T^{*}Q$, the cotangent lift, defined by $(g,(q,\dot q))\mapsto (T^{*}\psi_{g^{-1}}(q,\dot q))$. It enjoys the following properties \cite{Foundations}, \cite{manolo}:
\begin{itemize}
	\item $\psi^{T^{*}Q}$ is a symplectic action, meaning that $(\psi^{T^{*}Q}_g)^*\Omega=\Omega$, with $\Omega$ being the canonical symplectic $2$-form on $T^{*}Q$, $\Omega=dq\wedge dp$.
	\item It admits an $\hbox{Ad}^*$-equivariant momentum map 
	$
	J\colon T^*Q\to \mathfrak{g}^*$ given by
	$\langle J(q,p), \xi\rangle=\langle p, \xi_Q\rangle,\quad \forall \xi\in\mathfrak{g}$, where $\xi_Q(q)=d(\psi_{\exp(t\xi)}q)/dt$ is the infinitesimal generator of  $\xi\in \mathfrak{g}$, with $\mathfrak{g}$ the Lie algebra of $G$.
\end{itemize}
Likewise, $\psi^{TQ}$ denotes the tangent lift action on $TQ$, defined by $\psi^{TQ}_g=T\psi_{g}(q,\dot q)$.

To perform a hybrid reduction one needs to impose some compatibility conditions between the action and the hybrid system (see \cite{amesrouth} and \cite{amesham}). By an \emph{hybrid action} on the simple hybrid forced Lagrangian system $\mathscr{L}_F$  we mean a Lie group  action $\psi\colon G\times Q\to Q$ such that
\begin{itemize}
	\item $L$ is invariant under $\psi^{TQ}$, i.e. $L\circ \psi^{TQ}=L$.
	\item $\psi^{TQ}$ restricts to an action of $G$ on $S$.
	\item $\Delta$ is equivariant with respect to the previous action, namely $\Delta\circ \psi^{TQ}_g\mid_S=\psi^{TQ}_g\circ \Delta$.
\end{itemize}
Recall that $\psi^{TQ}$ admits an $\hbox{Ad}^*$-equivariant momentum map $J_L: TQ\to\mathfrak{g}^{*}$ given by $J_L=J\circ\mathbb{F} L$. This follows directly from the invariance of $L$, since it implies that $\mathbb{F} L$ is an equivariant diffeomorphism, i.e. $\mathbb{F} L\circ \psi^{TQ}_g=\psi_g\circ \mathbb{F}L$.

The hybrid equivalent of momentum map is the notion of \emph{hybrid momentum map} introduced in  \cite{amesrouth}, $J_L$ is an  \emph{hybrid momentum map} if the  diagram
\begin{equation}\label{diag1}
\begin{tikzcd}[column sep=1.5cm, row
    sep=1.2cm]
& \mathfrak{g}^* &\\
TQ \arrow[ur,"J_L"] & S \arrow[u,"J_L\mid_S"] \arrow[l,swap,hook',"i"] \arrow[r,"\Delta"] & TQ  \arrow[ul, swap,"J_L"]
\end{tikzcd} 
\end{equation}
commutes, where $i$ is the canonical inclusion from $S$ to $TQ$. 

We remind that (see \cite{manolo} for instance) by denoting $\{\phi^X_t\}$ the flow of a vector field $X$ on $Q$, we can also define the \textit{complete lift}  $X^c$ of $X$ in terms of its flow. We say that $X^c$ is the vector field on $TQ$ with flow $\{T\phi^X_t\}$. In other words, $X^c(v_q)=\left.\frac{d}{dt}\right|_{t=0}\left(T_q\phi^X_t(v_q)\right)$, or in local coordinates, $X^c=X^i\frac{\partial}{\partial q^i}+\dot{q}^j\frac{\partial X^i}{\partial q^j}\frac{\partial}{\partial \dot{q}^i}$.

For the Lagrangian side, one needs a further regularity condition, sometimes referred to as \emph{$G$-regularity}. Precisely, one has the following definition~\cite{Routhstages} (for an alternative, equivalent definition, see~\cite{quasi}). Let $L$ be an invariant Lagrangian on $TQ$ and denote by $\xi_Q$ the infinitesimal generator for the associated action. We say that $L$ is \emph{$G$-regular} if, for each $v_q\in TQ$, the map $\mathcal{J}_L^{v_q}:\mathfrak{g}\to \mathfrak{g}^*$, $\xi\mapsto J_L\left(v_q + \xi_Q(q)\right),\, v_q\in T_{q}Q$, is a diffeomorphism. In a nutshell, $G$-regularity amounts to regularity ``with respect to the group variables''. From now on we will assume that the Lagrangian is $G$-regular. In fact, this is always the case for mechanical Lagrangians.

Consider a simple hybrid forced Lagrangian system $\mathscr{L}_{F}=(TQ,X_{L}^{F}, S,\Delta)$ equipped with an hybrid action $\psi$ and $L$ invariant under $\psi^{TQ}$. We can apply a hybrid analog of the symplectic reduction Theorem for forced Lagrangian systems \cite{manololainzlopen} to the simple hybrid forced Lagrangian system $\mathscr{L}_{F}=(TQ,X_{L}^{F},S,\Delta)$ as follows: Consider the momentum map $J_L: TQ\to\mathfrak{g}^{*}$, given by $J_L(v_{q})(\xi)=\alpha_{L}(v_{q})(\xi_{Q}^{c}(v_{q}))$, where $\alpha_L=S^{*}(dL)$ being $S$ the vertical endomorphism on $TQ$ (see \cite{manolo}) locally given by $S=dq^{i}\otimes\frac{\partial L}{\partial \dot{q}^{i}}$ and denote by $\omega_{L}=-d\alpha_L$ the Poincar\'e-Cartan $2$-form \cite{manolo}, where the symbol $\otimes$ denotes a tensorial product. In addition, the invariance of $L$ implies the invariance of $\alpha_L$ and $\omega_L$ and the equivariance of the momentum map $J_L$ \cite{manolo}. 

For each $\xi\in\mathfrak{g}$ and $v_{q}\in TQ$, consider the function $J_L^{\xi}:TQ\rightarrow\R$ given by $J_L^{\xi}(v_{q})=\langle J_L(v_{q}),\xi\rangle.$ Let $\xi\in\mathfrak{g},$ then $J_L^{\xi}$ is a conserved quantity for $X_L^F$ if and only if $F(\xi_{Q}^{c})=0$ (see \cite{manololainzlopen}), where $\xi_{Q}^{c}$ denotes the complete lift of the vector field $\xi_Q$ given by the infinitesimal generator for the Lie group action $\psi$. In addition, the vector subspace of $\mathfrak{g}$ given by $\mathfrak{g}_{F}=\{\xi\in \mathfrak{g}:F(\xi_{Q}^{c})=0,\,\, i_{\xi_{Q}^{c}}dF=0\}$ is a Lie subalgebra of $\mathfrak{g}$ (see \cite{manololainzlopen}).  In particular, for each $\xi\in\mathfrak{g}_{F}$, $\xi_{Q}^{c}$ is a symmetry of the forced Lagrangian system given by $L$ and $F$.

Let $G_{F}\subset G$ be the Lie subgroup generated by $\mathfrak{g}_{F}$ and $J_{F}:TQ\rightarrow \mathfrak{g}_{F}^{*}$ the reduced hybrid momentum map. Let $\mu\in\mathfrak{g}_{F}^*$ be a \textit{hybrid regular value} of $J_{F},$ which means that $\mu$ is a regular value of both $J_{F}$ and $J_{F}\mid_S$ and let $(G_{F})_{\mu}$ be the isotropy subgroup in $\mu$.  Note that, since $L$ is invariant under $\psi^{TQ}$ and $G$-regular, then:
\begin{enumerate}[(i)]
\item  The reduced space $M_{\mu}\coloneqq J_{F}^{-1}(\mu)/(G_{F})_{\mu}$ is a symplectic manifold, with symplectic structure $\omega_{\mu},$ uniquely determined by  $\pi^{*}_{\mu}\omega_{\mu}=i^{*}_{\mu}\omega_{L},$ where $\pi_{\mu}:J_{F}^{-1}(\mu)\rightarrow M_{\mu}$ and $i_{\mu}:J_{F}^{-1}(\mu)\rightarrow TQ$ denotes the canonical projection and the  canonical inclusion, respectively. Moreover,  $J_{F}^{-1}(\mu)$ is a submanifold of $TQ$ and $X_{L}^F$ is tangent to it. 
 \item $L$ induces a function $R^\mu:M_{\mu}\rightarrow\R$ defined by $R^\mu\circ\pi_{\mu}=L\circ i_{\mu}.$
 \item $J_{F}\mid_{S}^{-1}(\mu)\subset S$ is $(G_{F})_{\mu}-$invariant and hence reduces to a submanifold of the reduced space which we denote $S_{\mu}\subset J_{F}^{-1}(\mu)/(G_{F})_{\mu}.$
 \item Again, using invariance, $\Delta$ reduces to a map $\Delta_{\mu}:S_{\mu} \rightarrow J_{F}^{-1}(\mu)/(G_{F})_{\mu}.$
 \item $F$ induces a reduced 1-form $F_{\mu}$ on $M_{\mu},$ uniquely determined by  $\pi^{*}_{\mu}F_{\mu}=i^{*}_{\mu}F$. 
\end{enumerate}

\vspace{.2cm}

A case of special interest with regards to applications is when $Q=\mathbb{S}^1\times M,$ where $M$ is called the shape space and the action is simply $(\theta,x)\mapsto(\theta+\alpha,x).$ This is often the situation when dealing with simple models of bipedal walkers, see e.g. \cite{ames2006}. From now on, we will assume we work in this setting. While this is indeed a strong assumption, it is always the case locally, so as long as it applies to the domain of interest of an specific problem the procedure applies. 

The forced Lagrangian system has a cyclic coordinate $\theta$, i.e., $L$ is a function of the form $L(\dot{\theta},x,\dot{x})$, and the forced $F$ is such that $F_{\theta}=0$ and $F_{x}$ is independent of $\theta$ with $F=F_{\theta}d\theta+F_{x}dx$. The conservation of the momentum map $J_{F}=\mu$ reads $\frac{\partial L}{\partial \dot{\theta}}=\mu,$ and one can use this relation to express $\dot{\theta}$ as a function of the remaining -non cyclic- coordinates and their velocities, and the prescribed regular value of the momentum map $\mu.$ We point out that it is at this stage that $G$-regularity of $L$ is used: it guarantees that $\dot{\theta}$ can be worked out in terms of $x$, $\dot{x}$ and $\mu$. If one chooses the cannonical flat connection on $Q\rightarrow Q/\mathbb{S}^1=M,$ then the Routhian can be computed as
\begin{equation} \label{5}
R^\mu(x,\dot x)=\left[L(\dot{\theta},x,\dot x)-\mu\dot{\theta} \right]\Big{|}_{\dot{\theta}=\dot{\theta}(x,\dot x,\mu)},    
\end{equation}
where the notation means that we have everywhere expressed $\dot{\theta}$ as a function of $(x,\dot x,\mu)$. Note that (\ref{5})
is the classical definition of the Routhian \cite{routhdissipative}. Let us first consider the case in which the momentum map is preserved in the collisions with the switching surface (elastic case). We then have:
\begin{proposition}\label{(3)} In the situation above:
\begin{enumerate}[(a)]
 \item Any solution of $\mathscr{L}_{F}=(TQ,X^{F}_{L},S,\Delta)$ with momentum $\mu$ projects onto a solution of $\mathscr{L}_ {F}^{\mu}=(T(Q/\mathbb{S}^1),X_{R_F}^{\mu},S_{\mu},\Delta_{\mu})$.
 \item Any solution of $\mathscr{L}_ {F}^{\mu}=(T(Q/\mathbb{S}^1),X_{R_F}^{\mu},S_{\mu},\Delta_{\mu})$ is the projection of a solution of $\mathscr{L}_{F}=(TQ,X^{F}_{L},S,\Delta)$ with momentum $\mu$.
\end{enumerate}
\end{proposition}

\vspace{.2cm}

Collisions with the switching surface will, in general, modify the value of the momentum map (nonelastic
case). Therefore, if $\mathcal{J}=\{I_{i}\}_{i\in \Lambda}$ is the hybrid interval, the Routhian has to be defined in each $I_i$ taking into account the value of the momentum $\mu_i$ after the collision at time $\tau_i$. Note that this also has influence in the way the reset map $\Delta$ is reduced. This will be clarified in the examples below.

Let us denote: (1) $\mu_i$ the momentum of the system in $I_i=[\tau_i,\tau_{i+1}]$, (2) $\Delta_{\mu_i}$ the reduction of $\Delta$, and (3) ${S}_{\mu_i}$ the reduction of ${S}$, so there is a sequence of reduced simple hybrid Routhian systems. The fact that the momentum will, in general, change with the collisions makes the reconstruction procedure more challenging. If one wishes, as usual, to use a reduced solution to reconstruct the original
dynamics, one needs to compute the reduced hybrid data after each collision. This means that once the reduced solution has been obtained between two collison events, say at $t=\tau_{n} $ and $t=\tau_{n+1},$ one
should reconstruct this solution to obtain the new momentum after the collision at $\tau_{n+1}$ and use this
new momentum to build a new reduced hybrid system whose solution should be obtained until the next collision event at $\tau_{n+2}$ and so on. As usual, the reconstruction procedure from the reduced hybrid flow to the hybrid flow involves an integration at each stage in the previous diagram of the cyclic variable using
the solution of the reduced simple hybrid forced Lagrangian system. Essentially, this accounts to imposing the momentum constraint on the reconstructed solution.

\vspace{.2cm}


\begin{example}
Continuing with Examples \ref{example} and \ref{ex2} if we square both sides of \eqref{eq12}, and noting that
		\[
		2x\dot x+2y\dot y=\frac{d}{dt}(x^2+y^2)= \frac{d}{dt}(r^2)=2r\dot r,
		\]
	we have $(\dot x^+)^2= (\dot x^-)^2+(2r\dot r^{-})^2 x^2-4x\dot x^{-}r\dot r^{-},$
		and symmetrically to $\dot y^{+}$. Add $(\dot x^+)^2+(\dot y^+)^2$. We can conclude
		\begin{align*}
		(\dot r^{+})^2&=(\dot r^{-})^2+(2r\dot r^{-})^2(x^2+y^2)-4(x\dot x^{-}+y\dot y^{-})r\dot r^{-}\\
		&=(\dot r^{-})^2+4r^2(\dot r^{-})(r^2-1).
		\end{align*}
		
		This means that, since the collision occurs at $r=1$, we have to $(\dot r^{+})^2=(\dot r^{-})^2$, then the solution that is obtained (physically) is $\dot r^{+}=-\dot r^{-}$. For $\theta=\arctan (y/x)$, we have 
		\begin{equation*}
		\dot \theta^+=\frac{1}{1+(y/x)^2} \left(\frac{\dot y^+ x-y \dot x^+ }{x^2}\right)
		=\frac{1}{r^2} \left(\dot y^- x - y \dot x^-\right)=\dot \theta^-.
		\end{equation*}
		where we have replaced the expression for $\dot x^+$, $\dot y^+$ and we used that $x^2+y^2=r^2$ and $\left(y \dot x^- -\dot y^- x\right)=-r^2\dot\theta^-$. It is understood that the ``minus'' square root is taken in $\dot{r}^{+}$ (the particle bounces on the boundary after the collision). The assumption of elastic collision implies, in particular, that the momentum map is preserved. This is clear since $r$ and $\dot\theta$ do not change with the collision. The Routhian, the reduced force and the reduced forced Euler-Lagrange equations are given in Example \ref{example}. The reduced reset map is determined by the expression for $\dot r^+$ (note that the expression drops to the quotient since it only involves $r$ and $\dot r$ ). The reduced switching surface is $S_{\mu}=\{r^2=1, \dot r>0\}$. One obtains the simple hybrid forced Routhian system $\mathscr{L}^{F}_{\beta}=(TQ_{\rm red},R^{\mu},S_\mu,\Delta_{\mu})$, with $Q_{\rm red}\simeq \R^+$ parametrized by the radial coordinate $r$.
\begin{figure}[h]
	\centering\includegraphics[height=3.8cm]{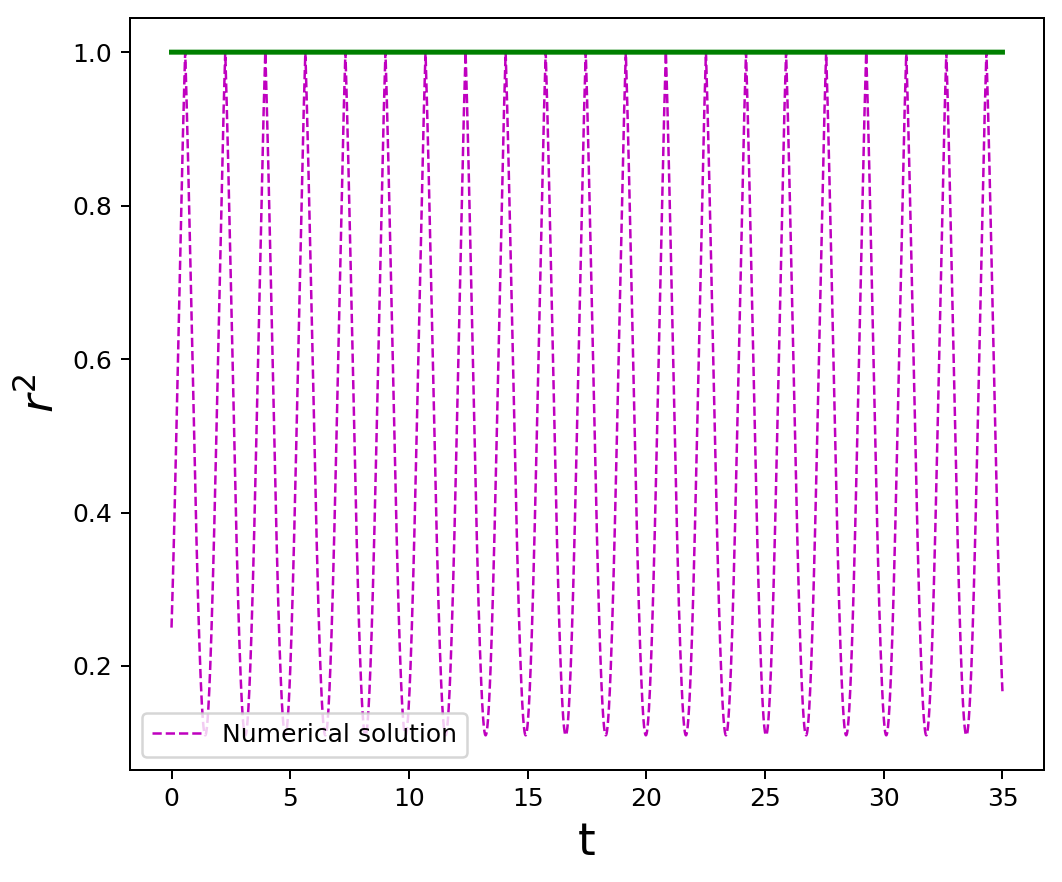}
	\centering\includegraphics[height=3.8cm]{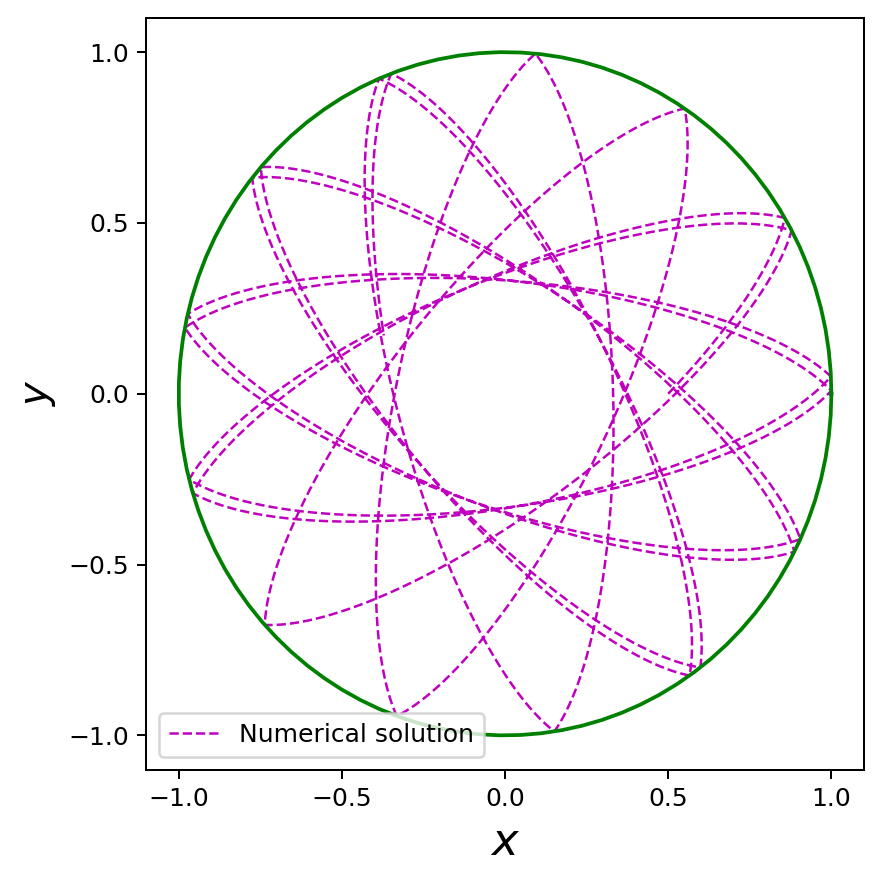}
	\caption{Simulation for $c=2$. The first figure  corresponds with the reduced trajectory while the second figure with the reconstructed solution.}\label{fig:2}
\end{figure}

Figures~\ref{fig:2} and~\ref{fig:3} show numerical results using \textsc{Python} for two different values of the dissipation parameter $c$. The remaining parameters are the same for both simulations: $m=1$, $r(0)=0.5$, $\dot r(0)=2$, $\theta(0)=0$ (rad) and $\dot \theta(0)=1$ (rad/s). The reduced dynamics is solved numerically (dashed purple line) and used to integrate (numerically) the reconstruction equation~$
\dot\theta=\frac{\mu} {m r^2}$, with $\mu$ determined from the initial conditions. Switching surfaces $S$ and $S_{\mu}$ are represented with a green solid line. Note also that the impact times on which the particle bounces are also obtained numerically. 
\begin{figure}[h]
	\centering\includegraphics[height=3.8cm]{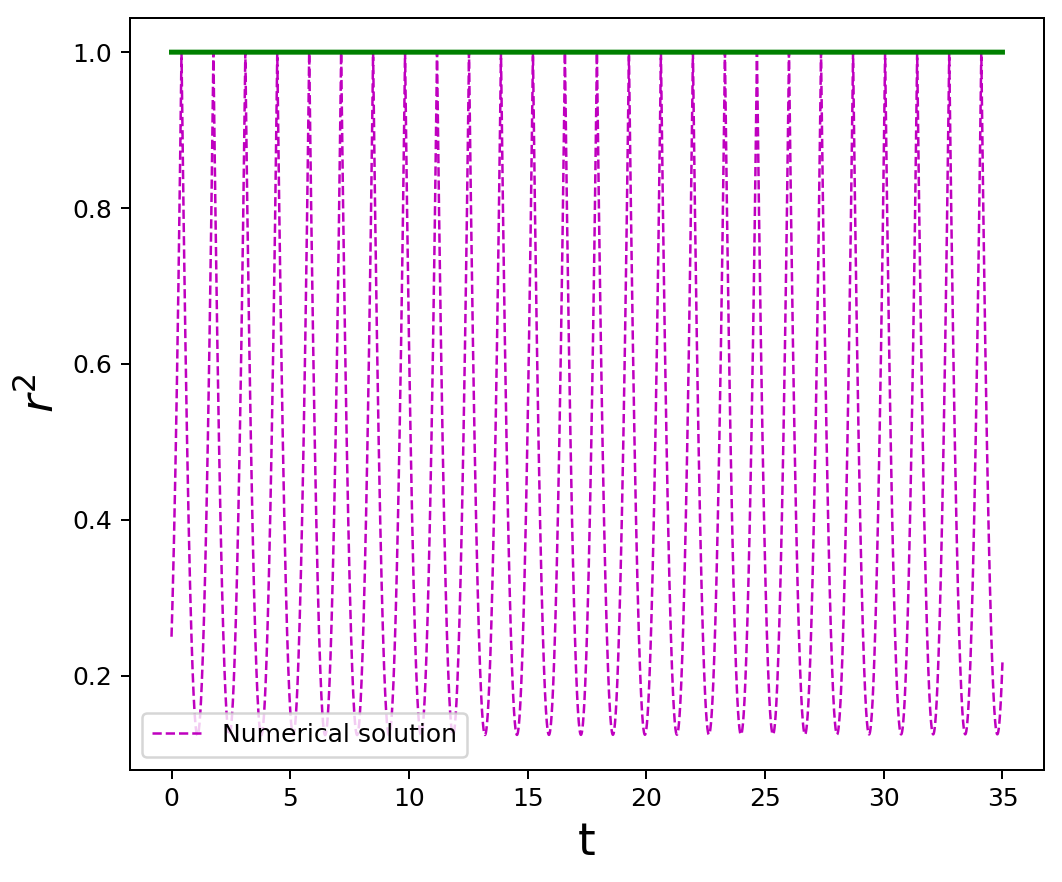}
	\centering\includegraphics[height=3.8cm]{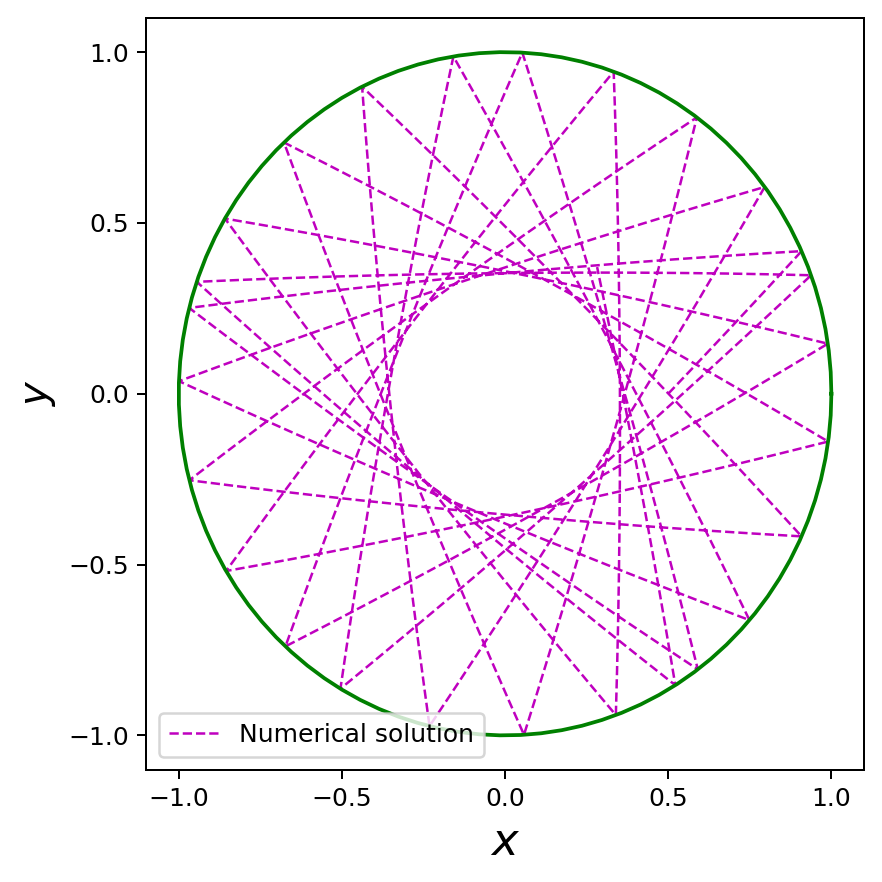}
	\caption{Simulation for $c=0.20$. The first figure corresponds with the reduced trajectory while the second figure with the reconstructed solution.}\label{fig:3}
\end{figure}

\end{example}

\section{Acknowledgment}
This work was supported by Ministerio de Ciencia e Innovación (Spain) under grant PID2019-106715GB-C21. The authors also acknowledge Manuela Gamonal for the help provided with the numerical simulations.
%
%
%
%
%

%
%


\begin{thebibliography}{10}

\bibitem{Foundations}
R.~Abraham and J.~E. Marsden.
\newblock {\em Foundations of mechanics}.
\newblock Benjamin/Cummings Publishing Co., Inc., Advanced Book Program,
Reading, Mass, 1978.

\bibitem{ames2006}
A.~Ames, R.~Gregg, E.~Wendel, and S.~Sastry.
\newblock \textit{On the geometric reduction of controlled three-dimensional bipedal
  robotic walkers}.
\newblock In 3rd Workshop on Lagrangian and Hamiltonian Methods for
  Nonlinear Control, 2006.



\bibitem{amesham} A. Ames, S. Sastry. \textit{Hybrid cotangent bundle reduction of simple hybrid mechanical systems with symmetry}. in Proceedings of the 25th American Control Conference Minneapolis MN 2006.
\bibitem{amesrouth} A. Ames, S. Sastry. \textit{Hybrid Routhian reduction of Lagrangian hybrid systems}. in Proceedings of the 25th American Control Conference Minneapolis MN 2006.

\bibitem{Bl}
A.~M. Bloch. \textit{Nonholonomic Mechanics and Control}, Series IAM. New York: Springer-Verlag, vol.~24. 2nd Edition 2015.

\bibitem{brogliato}  B.  Brogliato, \textit{Nonsmooth Impact Dynamics: Models, Dynamics and Control.} Springer-Verlag, 1996, vol. 220.
\bibitem{will} Clark, W. and Bloch, A. The Bouncing Penny and Nonholonomic Impacts. In Proceedings of the 2019 IEEE 58th Conference on Decision and Control (CDC) pp. 2114-2119, 2019

\bibitem{leo1}L. Colombo and D. V. Dimarogonas. Symmetry Reduction in Optimal Control of Multiagent Systems on Lie Groups. in IEEE Transactions on Automatic Control, vol. $65$, no. $11$, pp. $4973-4980$, 2020.
\bibitem{leo2} L. Colombo, F. Jimenez, and D. Mart\'in de Diego. Variational integrators for mechanical control systems with symmetries. Journal of Computational Dynamics, vol. 2, no. 2, pp. 193–225, 2015.

\bibitem{tdep} L. Colombo, M. E. Eyrea Iraz\'u, and Eduardo Garc\'ia-Torano Andr\'es. \textit{A note on Hybrid Routh reduction for time-dependent Lagrangian systems}. Journal of Geometric Mechanics 12(2), 309-321, 2020. arXiv: 2003.07484
\bibitem{NAHS} L. Colombo, M. E. Eyrea Iraz\'u. \textit{Symmetries and periodic orbits in simple hybrid Routhian systems}. Nonlinear Analysis: Hybrid Systems 36 (2020), 100857. arXiv:2001.08941
\bibitem{cortes}J. Cort\'es, M. de Le\'on, M. Mart\'in de Diego, S. Mart\'inez. \textit{Mechanical systems subjected to generalized non-holonomic constraints}.  R. Soc.  Lond.Proc.  Ser.  A Math.  Phys.  Eng.  Sci.  457(2007), 651-670, 2001.
\bibitem{cortes2}J. Cort\'es, A. Vinogradov.  \textit{Hamiltonian  theory  of  constrained impulsive motion}.  J. Math.  Phys.  47(4), 042905, 30 pp, 2006.

\bibitem{emmapoisson} M. E. Eyrea Irazú, L. Colombo, A. Bloch.  Reduction by Symmetries of Simple Hybrid Mechanical Systems.  Proceedings of the 7th IFAC Workshop on Lagrangian and Hamiltonian Methods in Nonlinear Control,  2021.
\bibitem{EmmaPhd}
E.~Eyrea~Irazu.
\newblock {\em Aspectos geométricos y numéricos de los sistemas mecánicos con términos magnéticos}.
\newblock PhD thesis. Facultad de Ciencias Exactas,
Universidad Nacional de La Plata, 2019.


\bibitem{Edu} E. Garc\'ia-Tora\~no, B. Langerock, F. Cantrijn. Aspects of reduction and transformation of Lagrangian systems with symmetry. J. Geom. Mech. 6 (2014), no. 1, 1-23.
\bibitem{hybridbook}  R. Goebel, R. Sanfelice, and A. Teel. Hybrid Dynamical Systems: modeling, stability, and robustness. Princeton University Press, 2012.

\bibitem{goldstein} H. Goldstein. \textit{Classical mechanics}. Addison-Wesley Publishing Co., second edition, 1980. Addison-Wesley Series in Physics.
\bibitem{grizzle2} J. Grizzle, S. Marcus. The structure of nonlinear control systems possessing symmetries. IEEE Transactions on Automatic Control, 30(3), 248-258, 1985.


\bibitem{grizzle} J.  Grizzle,  G.  Abba,  and  F.  Plestan,  \textit{Asymptotically  stable  walking for  biped  robots:  analysis  via  systems  with  impulse  effects}, IEEE Transactions on Automatic Control, vol. 46, no. 1, pp. 51-64, 2001.

	\bibitem{HoFuKoGu}  P. Holmes, R. Full, D. Koditschek, J.  
Guckenheimer. \textit{The dynamics of legged locomotion: models,  
	analyses, and challenges.} SIAM Review 48, no. 2, 207-304, 2006.
\bibitem{ibort} A. Ibort, M. de  Le\'on, E.  Lacomba, J.C.  Marrero, D. Mart\'in  de  Diego, P. Pitanga. \textit{Geometric  formulation  of Carnot's theorem}.  J. Phys.  A 34 (2001), no.  8, 1691-1712.
\bibitem{ibort2} A. Ibort, M. de  Le\'on, E.  Lacomba, D.  Mart\'in de  Diego, P. Pitanga. \textit{Mechanical  systems  subjected  to  impulsive  constraints}.  J. Phys.  A 30 (1997), no.  16, 5835-5854.

\bibitem{SHS} S. Johnson. \textit{Simple hybrid systems} Int. J. Bifurcation and Chaos, 04, 1655, 1994.

\bibitem{JK2}
 E.~Justh, P.~Krishnaprasad. Optimality, reduction and collective
  motion, Proc. R. Soc. A, 471 (2015), 20140606.

\bibitem{Leonard1}
 N.~Leonard, P.~Krishnaprasad. Motion control of drift-free,
  left-invariant systems on lie groups. IEEE Transactions on Automatic Control, 40 (1995), 1539--1554.
  \bibitem{manolo} M. de Le\'on and P. R. Rodrigues. Methods of differential geometry in analytical mechanics, ser. North-Holland Mathematics Studies 158. NorthHolland, 1989, isbn: 978-0-444-88017-8
\bibitem{manololainzlopen} de Le\'on, M., Lainz, M., $\&$ L\'opez-Gord\'on, A. (2021). Symmetries, constants of the motion, and reduction of mechanical systems with external forces. Journal of Mathematical Physics, 62(4), 042901.

	\bibitem{quasi}
	B.~Langerock, F.~Cantrijn, J.~Vankerschaver.
	\newblock Routhian reduction for quasi-invariant {L}agrangians.
	\newblock {\em J. Math. Phys.}, 51(2):022902, 20, 2010.
\bibitem{Routhstages}
B. Langerock, T. Mestdag, J. Vankerschaver.
\newblock Routh reduction by stages.
\newblock {\em SIGMA Symmetry Integrability Geom. Methods Appl.}, 7:Paper 109,
  31, 2011.


\bibitem{routhdissipative} L.A. Pars. A Treatise on Analytical Dynamics. Heinemann Educational Books,(1965).


\bibitem{bonnabel} A. Sarlette, S. Bonnabel, R. Sepulchre. Coordinated motion design on lie groups. IEEE Trans. Automatic Control, 55(5), 1047–1058, 2010.


\bibitem{tomlin}  C. Tomlin, Y. Ma and S. Sastry. Free flight in 2000: games on Lie groups, Proceedings of the 37th IEEE Conference on Decision and Control, vol.2, 2234-2239, 1998.



\bibitem{Biped-book}E. Westervelt, J. Grizzle, C. Chevallereau, J. Ho Choi,  and  B.  Morris. \textit{Feedback  control  of  dynamic  bipedal robot locomotion}.  Taylor $\&$ Francis/CRC, 2007.

\end{thebibliography}
\end{document}